\documentclass[12pt]{article}
\usepackage{amssymb}
\usepackage{amsmath}
\usepackage{graphics,url,color,epsfig}

\newcommand{\lap}{\mbox{$\bigtriangleup$}}

\newcommand{\ra}{{\mbox{$\rightarrow$}}}
\newcommand{\be}{\begin{equation}}
\newcommand{\ee}{\end{equation}}
\newtheorem{mrem}{Remark}
\newtheorem{mthm}{Theorem}

\newtheorem{mlem}{Lemma}

\newtheorem{lem}{Lemma}[section]

\begin{document}

\title{Monotonicity of solutions for fractional equations with De Giorgi type nonlinearities}

\author{Leyun Wu \quad Wenxiong Chen \thanks{ The corresponding author.} }

\date{\today}
\maketitle
\begin{abstract} In this paper, we develop a sliding method for the fractional Laplacian. We first obtain the key ingredients needed in the sliding method either in a bounded domain  or in the whole space, such as narrow region principles and maximum principles in unbounded domains.  Then using semi-linear equations involving the fractional Laplacian in both bounded domains and in the whole space, we illustrate how this new sliding method can be employed to obtain monotonicity of solutions. Some new ideas are introduced. Among which, one is to use Poisson integral representation of $s$-subharmonic functions in deriving the maximum principle, the other is to estimate the singular integrals defining the fractional Laplacians along a sequence of approximate maximum points by using a generalized average inequality. We believe that this new inequality will become a useful tool in analyzing fractional equations.
\end{abstract}
\bigskip

\bigskip

\section{Introduction}

The fractional Laplacian has attracted much attention recently. It has various applications in  anomalous diffusion and quasi-geostrophic
flows, turbulence and water waves, molecular dynamics, and relativistic quantum mechanics
of stars ( see \cite{BoG, CaV, Co, TZ} and the
references therein). It also has various applications in probability
and finance (see \cite{A, Be, CT}). In particular, the fractional Laplacian
can be understood as the infinitesimal generator of a stable L\'{e}vy diffusion process (see \cite{Be}). It  also has connections to conformal geometry, e.g. \cite{CG}.
\medskip

The fractional Laplacian in $\mathbb{R}^n$ is a nonlocal operator, taking the form
\begin{eqnarray}\label{eqFL}
(-\Delta)^s u(x) = C_{n, s} \, PV \int_{\mathbb{R}^n} \frac{u(x)-u(y)}{|x-y|^{n+2s}} dy,
\end{eqnarray}
 where $s$ is any real number between $0$ and $1$ and PV stands for the Cauchy principal value. It can be evaluated as
 \begin{eqnarray*}
  C_{n, s} \lim_{\epsilon \rightarrow 0^{+}}\int_{\mathbb{R}^{n}\backslash B_{\epsilon}(x)}
\frac{u(x)-u(y)}{|x-y|^{n+2s}}dy.
 \end{eqnarray*}
\medskip

In order that the integral on the right-hand side of (\ref{eqFL}) is well defined, we require that
$$ u \in C^{1,1}_{loc} \cap {\cal{L}}_{2s} := \{ u \in L^1_{loc} \mid \int_{\mathbb{R}^n} \frac{|u(x)|}{1 + |x|^{n+2s}} d x < \infty \}.$$
\medskip

 The non-locality of the fractional Laplacian makes it difficult to investigate.  To circumvent this, Caffarelli and Silvestre [CS] introduced  the extension method that reduced this nonlocal problem into a local one  in higher dimensions.
Another approach is the integral equations methods \cite{CLO}. After establishing the equivalence between a fractional equation and its corresponding integral equation, one can apply
the method of moving planes in integral forms or regularity lifting  to obtain the symmetry and regularity of solutions to the fractional equations.

These methods have been applied successfully to study equations involving the fractional Laplacian, and a series of fruitful results have been  obtained (see \cite{BCPS}, \cite{CZ}, \cite{CFY},  \cite{CLO1}, \cite{ZCCY}, and the references therein).

However, when applying the above two methods, sometimes one need to impose extra conditions on the solutions, which would not be necessary when one considers the pseudo-differential equation directly.
Moreover, they do not work for nonlinear nonlocal operators, such as the fractional \emph{p}-Laplacians (see \cite{CQ} for more details).

Hence it is more desirable to develop direct methods without going through extensions or integral equations. Direct methods of moving planes for the fractional Laplacian \cite{CLL} \cite {CLM}
and for fractional \emph{p}-Laplacians \cite{CL}
have been introduced, and have been applied to obtain symmetry, monotonicity, and non-existence of solutions for various semi-linear equations involving these nonlocal operators.
\medskip

In this paper, we introduce a direct sliding method for the fractional Laplacian.
The sliding method was developed  by Berestycki and  Nirenberg (\cite{BN1}-\cite{BN3}).
 It was used to establish qualitative properties of solutions for partial differential equations (mainly involving the regular Laplacian) such as symmetry, monotonicity, and uniqueness etc...
 The essential ingredients are different forms of maximum principles.  The main idea lies in comparing values of the solution of the equation at two different points, between which one point is obtained from the other by sliding the domain in a given direction, and then the domain is slid back to a critical position. While in the method of moving planes, one point is the reflection of the other.
\medskip

The following are some typical applications of the sliding method.
\medskip

In \cite{BN1}, Berestycki and  Nirenberg studied the monotonicity and uniqueness of the equation:
$$
\Delta u+ f(x,u, \nabla u)=0
$$
in a finite cylinder.

In \cite{BN2} and \cite{BCN1}, Berestycki, Caffarelli, and  Nirenberg studied the monotonicity, symmetry, and uniqueness of the equation:
$$
\Delta u+ \beta(x')u_n+f(x', u)=0
$$
in an infinite cylinder.
\medskip

In \cite{BN3} and  \cite{BCN2}, Berestycki, Caffarelli, and  Nirenberg investigated the monotonicity and uniqueness of the equation:
$$
\Delta u+ f(u)=0
$$
in a bounded domain, and generalized the results to fully nonliner equations.
\medskip

In \cite{DSV}, Dipierro, Soave and Valdinoci studied the over-determined problems of the type

  \begin{eqnarray*}
  \left\{\begin{array}{ll}
  (-\Delta)^s u =f(u) &  \mbox{ in } \Omega,\\
  u>0& \mbox{ in } \Omega,\\
  u=0 & \mbox{ in } \mathbb{R}^n\backslash \Omega,\\
({\partial}_\nu)_s u=const. & \mbox{ on } \partial\Omega.
\end{array} \right.
\end{eqnarray*}
where $\Omega$ is the region above the graph of a continuous function $\varphi:\mathbb{R}^{n-1}\rightarrow\mathbb{R}$:
\begin{eqnarray*}
\Omega:=\{x=(x',x_n)\in\mathbb{R}^n\mid x_n>\varphi(x')\} \ \ \text{with}\ \ x'=(x_1,\cdots,x_{n-1})\in \mathbb{R}^{n-1}.
\end{eqnarray*}
They first obtained the monotonicity of solutions and then showed that the $\Omega$ must be a half space under some appropriate conditions on $f$ and on $\Omega$.
\medskip

In this paper, we consider the fractional semi-linear equation
$$(-\Delta)^s u(x)=f(u(x))$$
in two different type of regions, bounded domains and the whole space.
\medskip

Before stating our main results, we introduce some notation. For
$$x = (x', {x_n}) \mbox{ with } x' = ({x_1},...,{x_{n-1}}) \in {\mathbb{R}^{n - 1}}$$ and  $\tau  \in \mathbb{R},$ let
$$u^\tau(x) =u(x', x_n+\tau)$$
and
$$w^\tau(x)=u^\tau(x)-u(x).$$
\medskip

Similar to the method of moving planes,  the narrow region principle is a key ingredient in the sliding method and it provides a starting position to slide the domain.
Hence, in this paper, we first establish
\medskip

\begin{mthm} [Narrow region principle]\label{mthm1}
 Let $D$ be a bounded narrow region in $\mathbb{R}^n.$  Suppose that  $u \in {{\cal{L}}_{2s}(\mathbb{R}^n)} \cap C_{loc}^{1,1}(D ),$ $w^\tau$ is lower semi-continuous on $\overline{D},$  and satisfies
 \begin{eqnarray}\label{eqF}\left\{\begin{array}{ll}
 ( - \Delta )^s{w^\tau }(x) + c(x){w^\tau }(x) \geqslant 0 &\quad \mbox{ in } D, \\
{w^\tau }(x)\geq 0 &\quad \mbox{ in } D^c,
\end{array} \right. \end{eqnarray}
with $c(x)$ bounded from below in $D$.

Let $d_n(D)$ be the width of $D$ in the $x_n$ direction, in which we assume that $D$ is narrow:
 \begin{eqnarray}\label{Fwidth}
d_n(D)|\inf_D c(x)|^{\frac{1}{2s}}\leq C,
 \end{eqnarray}

Then
\begin{eqnarray}\label{Feq1.4}
 {w^\tau }(x) \geqslant 0 \mbox{ in } D.
\end{eqnarray}
More strongly, we have
\begin{eqnarray}\label{Feq1.4-1}
 \mbox{ either } {w^\tau}(x) > 0 \mbox{ in } D \, \mbox{ or } w^\tau (x) \equiv 0 \mbox{ in } \mathbb{R}^n. \end{eqnarray}
\end{mthm}
\medskip

For a bounded region $\Omega$ which is convex in $x_n$-direction, let
$$\Omega^\tau=\Omega-\tau e_n \mbox{ with } e_n =(0, \cdots, 0, 1),$$
which is obtained by sliding $\Omega$ downward $\tau$ units.
\medskip

It is obvious that when $\tau$ is sufficiently close to the width of $\Omega$ in $x_n$-direction,
$\Omega \cap \Omega^\tau$ is a narrow region. Then under some monotonicity conditions on the solution $u$ in the complement of $\Omega$, we will be able to apply Theorem \ref{mthm1} to conclude that
\begin{eqnarray}\label{wtau}
w^\tau(x)\geq 0,\;\; x \in \Omega \cap \Omega^\tau.
\end{eqnarray}
\medskip

This provides a starting position to slide the domain $\Omega^\tau$. Then in the second step, we slide $\Omega^\tau$ back upward as long as inequality (\ref{wtau}) holds to its limiting  position. If we can slide the domain all the way to $\tau =0$, then we conclude that the solution is monotone increasing in
$x_n$-direction.
\medskip

To ensure the two steps, we need to impose the exterior condition on $u.$ Let
 \begin{eqnarray}\label{eqfnrqb}
u(x)=\varphi(x) \;\mbox{  in  } \;\Omega^c,
\end{eqnarray}
and assume that
\medskip

{\textbf{(H):}}
For any three points $x=(x', x_n), y=(x', y_n)$ and $z=(x', z_n)$ lying on a segment parallel to the $x_n$-axis, $y_n<x_n<z_n,$ with $y, z \in \Omega^c$, we have
\begin{eqnarray}\label{eq1.7}
\varphi(y)<u(x)<\varphi(z),\;\; \mbox{ if } x\in \Omega
\end{eqnarray}
and
\begin{eqnarray}\label{eq1.8}
\varphi(y)\leq \varphi(x)\leq \varphi(z), \;\;\mbox{ if } x\in \Omega^c.
\end{eqnarray}
\medskip

\begin{mrem}
The same monotonicity conditions (\ref{eq1.7}) and (\ref{eq1.8}) (with $\Omega^c$ replaced by $\partial \Omega$) were assumed in \cite{BN3}.
\end{mrem}

By employing the sliding method, we obtain the monotonicity of solutions for fractional equations in bounded domains.
\medskip

\begin{mthm}\label{mthm2}
Let $\Omega$ be a bounded domain of $\mathbb{R}^n$ which is convex in $x_n$-direction. Assume that $u \in {\cal{L}}_{2s}(\mathbb{R}^n)\cap C_{loc}^{1,1}(\Omega)$ is a solution of
\begin{eqnarray}\label{eqfnrp}
\left\{\begin{array}{ll}
(-\Delta)^s u(x)=f(u(x))  &\mbox{ in } \Omega, \\
u(x)=\varphi(x)  & \mbox{ on } \Omega^c,
\end{array}
\right.\end{eqnarray}
and satisfies (H). The function $f$ is supposed to be Lipschitz continuous.
Then $u$ is monotone increasing with respect to $x_n$ in $\Omega,$ i.e. for any  $\tau >0$,
$$
u(x', x_n+\tau) >u(x', x_n) \;\mbox{ for all}\; (x', x_n), (x', x_n+\tau) \in \Omega.
$$
\end{mthm}
\medskip

To apply the sliding method on unbounded domains, a maximum principle plays an important role.
We prove
\medskip

\begin{mthm}\label{mthm3} ({\em Maximum principle in unbounded domains})
Let $D$ be an open set in $\mathbb{R}^n,$ possibly unbounded and disconnected, suppose that
 \begin{eqnarray}\label{domain}
 \mathop{\underline{lim}}\limits_{k\rightarrow \infty}\frac{|D^c\cap (B_{2^{k+1}}(q)\backslash B_{2^k}(q))|}{|B_{2^{k+1}}(q)\backslash B_{2^k}(q)|}>0,
 \end{eqnarray}
 where $q$ is any point in $D$. Let $u \in {{\cal L}_{2s}(\mathbb{R}^n)} \cap C_{loc}^{1,1}(D)$ be bounded from above, and satisfies
\begin{eqnarray}\label{eqFIM}\left\{\begin{array}{ll}
 ( - \Delta )^s u+c(x) u(x) \leq 0, & \mbox{ at the points in } D \mbox{ where } u(x)>0, \\
 u(x) \leq 0,& \mbox{ in } \mathbb{R}^n \backslash D
\end{array} \right. \end{eqnarray}
for some nonnegative function $c(x)$. Then
$$u(x) \leq 0 \; \mbox{ in }  \;D.$$
\end{mthm}
\begin{mrem} \label{eg} 
Dipierro, Soave and Valdinoci in \cite{DSV} proved this theorem by using Silvestre's growth lemma (\cite{S}) under the {\em exterior cone condition} that the complement of $D$ contains an infinite open connected cone.

Here we introduce a new idea in the proof--using the Poisson representation of $s$-subharmonic functions,
and thus significantly weakens the {\em exterior cone condition} to condition (\ref{domain}). To illustrate this, we list some typical examples of $D$ which satisfy our condition (\ref{domain}), but obviously does not satisfy the {\em exterior cone condition}:
\medskip

(1) Stripes: $D=\left\{ x \mid 2k< x_n<2k+1, k=0, \pm 1, \pm 2, \cdots\right\}(disconnected).$
\medskip

(2) Annulus: $D=\left\{ x \mid 2k< |x|<2k+1, k=0, 1, 2, \cdots\right\}(disconnected).$
\medskip

(3) Archimedean spiral: See the figure (connected).

\begin{figure}[htbp] \begin{minipage}[t]{0.3\textwidth}
 \centering \includegraphics[width=\textwidth]{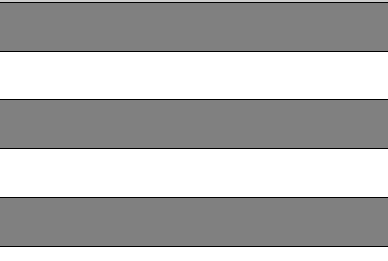} \caption{Stripes}

 \end{minipage}
 \begin{minipage}[t]{0.25\textwidth}
 \centering \includegraphics[width=\textwidth]{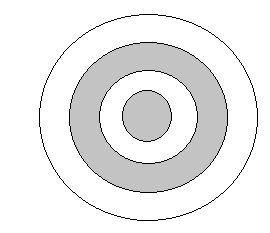}
 \caption{Annulus}
  \end{minipage}
   \begin{minipage}[t]{0.32\textwidth}
 \centering \includegraphics[width=\textwidth]{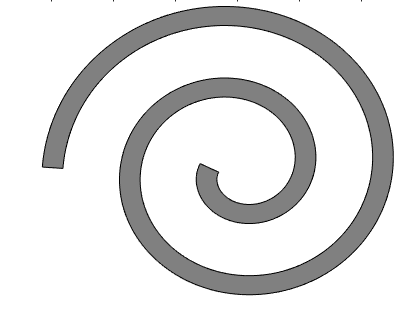}
 \caption{Archimedean spiral}
  \end{minipage}
  \end{figure}

\end{mrem}

Here $D$ are the shaded regions.
\medskip

\bigskip

\bigskip

\bigskip

\vspace{3in}

Next, as an application of the sliding method, we derive the monotonicity of solutions for fractional semi-linear equations in the whole space.
\medskip

\begin{mthm}\label{mthm4}
Let $u \in {\cal L}_{2s}(\mathbb{R}^n)\cap C_{loc}^{1,1}(\mathbb{R}^n)$ be a solution of
\begin{eqnarray}\label{eq1.1}
(-\Delta)^s u(x) =f(u(x)),\; \; x \in  \mathbb{R}^n,
\end{eqnarray}
$$|u(x)|\leq 1,$$ and
\begin{eqnarray}\label{eq1.2}
u(x', x_n)  \mathop {\longrightarrow }\limits_{{x_n\rightarrow \pm \infty}}\pm 1 \; \mbox{ uniformly in }\; x'=(x_1, \cdots, x_{n-1}).
\end{eqnarray}
Assume that $f(\cdot)$ is continuous in $[-1, 1]$ and there exists $\delta >0$ such that
\begin{eqnarray}\label{eq1.3}
f \;\mbox{ is nonincreasing on }\;[-1, -1+\delta] \;\mbox{ and on } \;[1-\delta, 1].
\end{eqnarray}

Then $u(x)$ is increasing with respect to $x_n,$ and furthermore, it depends on $x_n$ only.
\end{mthm}

\begin{mrem}
Theorem \ref{mthm4} is closely related to the De Giorgi conjecture (\cite{DG}), which states that:
\medskip

If $u$ is a solution of
\be -\lap u(x) = u(x)-u^3(x), \;\;\; x \in \mathbb{R}^n
\label{DG1}
\ee
such that
$$|u(x)| \leq 1, \; \lim_{x_n \ra \pm \infty} u(x',x_n) = \pm 1 \; \mbox{ for all }\; x' \in \mathbb{R}^{n-1}, \; \mbox{ and }\; \frac{\partial u}{\partial x_n} >0.$$
Then there exists a vector ${\bf a} \in \mathbb{R}^{n-1}$ and a function $u_1 : \mathbb{R} \ra \mathbb{R}$ such that
$$ u(x', x_n) = u_1({\bf a} \cdot x' + x_n), \;\; \forall \, x \in \mathbb{R}^n.$$

{\bf Note} the model function $f(u)=u-u^3$ does satisfy condition (\ref{eq1.3}).
\end{mrem}

In \cite{BHM}, Berestycki, Hamel, and Monneau proved the same monotonicity result for the following equation
$$
-\Delta u=f(u) \;\mbox{ in }\; \mathbb{R}^n
$$
under the same conditions as in Theorem \ref{mthm4} and additional conditions that
\begin{equation} \label{a10}
f=f(u) \mbox{ is Lipschitz continuous in } [-1, 1] \mbox{ and } f(\pm 1)=0.
\end{equation}
These conditions on $f$ were also required in \cite{DSV} when they considered fractional equations.

In this paper, we apply a new idea.
Instead of using the traditional approach--estimating along a sequence of equations in the whole domain $\Omega$ as in \cite{BHM} and \cite{DSV}, we estimate the singular integrals defining
$$(-\Delta)^s u (x) - (-\Delta)^s u^\tau (x)$$ along a sequence of approximate maximum points. This way, we will be able to weaken conditions
(\ref{a10}) and only need to assume that $f$ is continuous in Theorem \ref{mthm4}. In this process, the following average inequality plays a key role.

\begin{mlem}[A generalize average inequality]
Suppose that  $u \in {{\cal{L}}_{2s}(\mathbb{R}^n)} \cap C_{loc}^{1,1}(\mathbb{R}^n ),$ and $\bar{x}$ is a maximum point of $u$ in $\mathbb{R}^n.$ Then for any $r>0,$ we have
\begin{eqnarray*}
\frac{C_0}{C_{n, s}}r^{2s}(-\Delta)^{s}u(\bar{x})+ C_0 \int_{B_r^c(\bar{x})}\frac{r^{2s}}{|\bar{x}-y|^{n+2s}}u(y)dy\geq u(\bar{x}),
\end{eqnarray*}
where $C_0$ satisfies
$$
 C_0 \int_{B_r^c(\bar{x})}\frac{r^{2s}}{|\bar{x}-y|^{n+2s}}dy=1.
$$
\end{mlem}

In the special case when $u$ is $s$-subharmonic at point $\bar{x}$, the above inequality becomes
\begin{equation} u(\bar{x}) \leq \int_{B_r^c(\bar{x})} u(y) \, d \mu(y)
\label{AB1}
\end{equation}
with
$$ \int_{B_r^c(\bar{x})} d \mu(y) =1.$$
Here the integral on the right hand side of (\ref{AB1}) is actually a weighted average value of $u$ outside the ball
$B_r(\bar{x})$. This inequality can be conveniently used to prove maximum principles in unbounded domains such as Theorem \ref{mthm3}.
And we believe that it will become an effective tool in analyzing fractional equations.

\section{Narrow Region Principle and Monotonicity}

In this section, we prove Theorem \ref{mthm1} and Theorem \ref{mthm2}.
\medskip

{\textbf{Proof of Theorem \ref{mthm1}.}
\smallskip

Suppose (\ref{Feq1.4}) is not valid, then the lower semi-continuity of $w^\tau$ in $\bar{D}$ implies that there exists a point $x_0$ such that
 \begin{eqnarray}\label{Feq3.1}
 w^\tau(x_0)=\min_{\bar{D}} w^\tau(x) <0.
 \end{eqnarray}
By (\ref{eqF}) and (\ref{Fwidth}), we have
\begin{eqnarray}\label{Feq3.2}
 &&(-\Delta)^s w^\tau(x_0)+c(x_0)w^\tau(x_0)\nonumber \\
 &=&C_{n, s}PV \int_{\mathbb{R}^n}\frac{w^\tau(x_0)-w^\tau(y)}{|x_0-y|^{n+2s}}dy+c(x_0)w^\tau(x_0)\nonumber\\
 &\leq& C_{n, s}w^\tau(x_0) PV \int_{D^c}\frac{1}{|x_0-y|^{n+2s}}dy +\inf_D c(x) \; w^\tau(x_0)\nonumber\\
 &\leq& w^\tau(x_0)\left(\frac{C}{d_n^{2s}(D)} -\inf_D c(x)\right) \nonumber\\
   &<&0,
\end{eqnarray}
 where  $d_n(D)$ denotes the width of $D$ in the $x_n$ direction, and the second inequality from the bottom holds due to an argument in \cite{CLL}.

 Inequality (\ref{Feq3.2}) contradicts (\ref{eqF}),  and  we thus conclude that
$$w^\tau(x) \geq 0,  x \in D.$$

Based on this result, if $w^\tau(x)=0$ at some point $x\in D,$ then $x$ is a minimum point of $w^\tau$ in $D.$ If  $w^\tau \not\equiv 0$ in $\mathbb{R}^n,$ then we have
$$
(-\Delta)^s {w^\tau}(x)=C_{n,s}PV \int_{\mathbb{R}^n}\frac{{w^\tau}(x)-{w^\tau}(y)}{|x-y|^{n+2s}}dy<0.
$$
This contradicts
$$
 (-\Delta)^s {w^\tau}(x)= (-\Delta)^s {w^\tau}(x)+c(x){w^\tau}(x)\geq 0.
$$
 Therefore, we have
$$\mbox{ either } {w^\tau}(x) > 0 \mbox{ in } D \, \mbox{ or } w^\tau (x) \equiv 0 \mbox{ in } \mathbb{R}^n.$$

 This completes the proof of Theorem \ref{mthm1}.
\medskip

 Applying Theorem \ref{mthm1}, we derive the following monotonicity result of Theorem \ref{mthm2}.
\medskip

To better illustrate the idea of the sliding method, we only present the proofs for $D$ when it is an ellipsoid or a rectangle. When $D$ is an arbitrary bounded domain of $\mathbb{R}^n$ which is convex in $x_n$-direction, the proof is entirely similar.
\medskip

{\textbf{Proof of Theorem \ref{mthm2}.}
\smallskip

For $\tau\geq 0,$  denote
 $$u^\tau(x)=u(x', x_n+\tau).$$
 It is defined on the set $\Omega^\tau=\Omega-\tau e_n$ which is obtained from $\Omega$ by sliding it downward a distance $\tau$ parallel to the $x_n$-axis, where $e_n=(0, \cdots, 0, 1).$ Set
$$D^\tau:=\Omega^\tau \cap \Omega,$$
$$\tilde{\tau}=\sup\{\tau \mid \tau>0, D^\tau\neq \emptyset\},$$
 and
$$
w^\tau(x)=u^\tau(x)-u(x),\;\; x \in D^\tau.
$$
$u^\tau$ satisfies the same equation (\ref{eqfnrp}) in $\Omega^\tau$ as $u$ does in $\Omega,$ then $w^\tau$ satisfies
\begin{eqnarray}\label{eq3.12-0}
 ( - \Delta )^s{w^\tau(x) } = c^\tau(x){w^\tau }(x) \;\mbox{ in } D^\tau,
 \end{eqnarray}
where $$c^\tau(x)=\frac{f(u^\tau(x))-f(u(x))}{u^\tau(x)-u(x)}$$ is some $L^\infty$ function satisfying
$$c^\tau(x)\leq C, \;\;\forall x \in D^\tau.$$
\medskip

 The main part of the proof consists in showing that
 \begin{eqnarray}\label{fconclusion}
 w^\tau (x)>0, \;\; x \in  D^\tau, \;\mbox{ for any }\; 0<\tau<\tilde{\tau},
 \end{eqnarray}
this means precisely that $u$ is strictly increasing in the $x_n$ direction.

\begin{figure}[!ht]
\begin{center}
\begin{minipage}{0.40\textwidth}\centering
\epsfig{figure=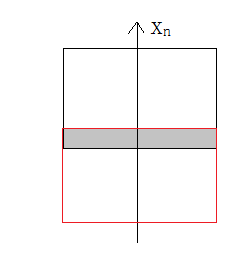,width=5.1cm}
\end{minipage}
\end{center}
\caption{The black $=\partial\Omega,$ the red$=\partial \Omega^\tau$, and the shaded region$=narrow \;region$.
\label{picture4}}
\end{figure}

\medskip

\emph{Step 1.}
\medskip

We show that
\begin{eqnarray}\label{start position}
 w^\tau(x) \geq 0 \;\mbox{ for } \tau \mbox{ sufficiently close to } \tilde{\tau} \mbox{ when } D^\tau \mbox{ is narrow}.
\end{eqnarray}
This can be derived directly from Theorem \ref{mthm1}.

\begin{figure}[!ht]
\begin{center}
\begin{minipage}{0.40\textwidth}\centering
\epsfig{figure=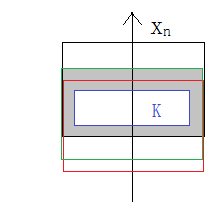,width=5.3cm}
\end{minipage}
\end{center}
\caption{The black $=\partial\Omega,$ the red$=\partial\Omega^\tau$, the green$=\partial\Omega^{\tau-\varepsilon},$  the blue $=\partial K$ and the shaded region$=narrow\;region$ .\label{picture5}}
\end{figure}

\medskip
\emph{Step 2.}
\medskip

 Inequality (\ref{start position}) provides a starting point, from which we can carry out the sliding.  Now we decrease $\tau$ as long as inequality (\ref{start position}) holds to its limiting position. Define
 $$
 \tau_0=\inf \left\{\tau \mid w^{\tau}(x)\geq 0, \; x \in D^\tau;\; 0<\tau<\tilde{\tau}\right\}.
 $$
 We prove that
 $$
 \tau_0=0.
 $$

 Otherwise, suppose that $\tau_0>0,$ we show that the domain $\Omega^\tau$ can be slid upward a little bit more and we still have
 \begin{eqnarray}\label{ww}
w^\tau (x)\geq 0, \;\; x \in  D^\tau, \;\mbox{ for any }\; \tau_0-\varepsilon<\tau \leq \tau_0,
\end{eqnarray}
which contradicts the definition of $\tau_0.$

 Since
 $$w^{\tau_0}(x)\geq 0, \;\;x \in D^{\tau_0},$$
 and
 $$w^{\tau_0}(x)>0, \;\; x \in \Omega \cap \partial D^{\tau_0}.$$
then
  $$w^{\tau_0}(x) \not\equiv 0, \;\; x\in D^{\tau_0}.$$

If there exists a point $x$ such that $w^{\tau_0}(x)= 0$, then $x$ is the minimum point and
$$
(-\lap)^s w^{\tau_0}(x)=C_{n, s}PV \int_{\mathbb{R}^n}\frac{w^{\tau_0}(x)-w^{\tau_0}(y)}{|x-y|^{n+2s}}dy<0.
$$
This contradicts $$
(-\lap)^s w^{\tau_0}(x)=f(u^{\tau_0}(x))-f(u(x))=0.
$$
Therefore,
\begin{eqnarray}\label{0ww}
w^{\tau_0}(x)>0, \;\; x\in D^{\tau_0}.
\end{eqnarray}

Now we can carve out of $D^{\tau_0}$ a closed set $K\subset D^{\tau_0}$ such that $D^{\tau_0}\backslash K$ is narrow. By (\ref{0ww}),
$$
w^{\tau_0}(x)\geq C_0>0 \mbox{ in } K.
$$
From the continuity of $w^\tau$ in terms of $\tau,$ we have for small $\varepsilon>0,$
$$w^{{\tau_0}-\varepsilon}(x) \geq 0 \mbox{ in } K.$$
In addition, we obtain from (H) that
$$
w^{{\tau_0}-\varepsilon}(x) \geq 0 \mbox{ in } (D^{\tau_0-\varepsilon})^c.
$$
Since $(D^{{\tau_0}-\varepsilon}\backslash K)^c=K\cup (D^{\tau_0-\varepsilon})^c,$  then
\begin{eqnarray}\label{eqww}
\left\{
\begin{array}{ll}
( - \Delta )^s{w^{{\tau_0}-\varepsilon} }(x) - c^{{\tau_0}-\varepsilon}(x){w^{{\tau_0}-\varepsilon} }(x) = 0, & x \in D^{{\tau_0}-\varepsilon}\backslash K,\\
{w^{{\tau_0}-\varepsilon } }(x)\geq 0, & x \in  (D^{{\tau_0}-\varepsilon}\backslash K)^c.
\end{array}
\right. \end{eqnarray}
It then follows from Theorem \ref{mthm1} that (\ref{ww}) is valid.
Therefore, we have reached a contradiction and we have
\begin{eqnarray}\label{ww0}
 w^\tau (x)\geq 0, \;\; x \in  D^\tau, \;\mbox{ for any }\; 0<\tau<\tilde{\tau},
 \end{eqnarray}

 Since
  $$w^{\tau}(x) \not\equiv 0, \;\; x\in D^{\tau}, \;\mbox{ for any }\; 0<\tau<\tilde{\tau},$$
if there exists a point $x^o$ such that $w^{\tau}(x^o)= 0$, then $x^o$ is the minimum point and
$$
(-\lap)^s w^{\tau}(x^o)=C_{n, s}PV \int_{\mathbb{R}^n}\frac{w^{\tau}(x^o)-w^{\tau}(y)}{|x^o-y|^{n+2s}}dy<0.
$$
This contradicts $$
(-\lap)^s w^{\tau}(x^o)=f(u^{\tau}(x^o))-f(u(x^o))=0.
$$
Hence, we
 arrived at (\ref{fconclusion}).

This completes the proof of Theorem \ref{mthm2}.
\medskip

\section{Maximum principle in unbounded domains and applications}

\subsection{The proof of Theorem \ref{mthm3}}

Denote
$$u^+:=\max\{u, 0\}, \;\; u^-:=\min\{u, 0\}$$
and
$$D^+:=\{x \in D : u(x) >0\}$$.

We prove it in three steps.

\emph{Step 1.}
\medskip

 In this step, we show that
\begin{eqnarray}\label{eq2.1}
(-\Delta)^s u^+ \leq 0 \;\mbox{ in }\; D^+.
\end{eqnarray}
In fact, if $x \in D^+,$ then
\begin{eqnarray*}
 &&(-\lap)^s u^+ (x)\\
 &=&C_{n, s}PV \int_ {\mathbb{R}^n} \frac{u^+(x)-u^+(y)}{|x-y|^{n+2s}}dy\\
 &=&C_{n, s}PV \int_ {D^+} \frac{u^+(x)-u^+(y)}{|x-y|^{n+2s}}dy+C_{n, s}PV \int_ {\mathbb{R}^n \backslash D^+} \frac{u^+(x)-u^+(y)}{|x-y|^{n+2s}}dy\\
  &=&C_{n, s}PV \int_ {D^+} \frac{u(x)-u(y)}{|x-y|^{n+2s}}dy+C_{n, s}PV \int_ {\mathbb{R}^n \backslash D^+} \frac{u(x)-u(y)}{|x-y|^{n+2s}}dy\\
  &&+C_{n, s}PV \int_ {\mathbb{R}^n \backslash D^+} \frac{u^-(y)}{|x-y|^{n+2s}}dy\\
  &=&(-\Delta)^s u(x)+C_{n, s}PV \int_ {\mathbb{R}^n \backslash D^+} \frac{u^-(y)}{|x-y|^{n+2s}}dy\\
  &\leq&(-\Delta)^s u(x).
\end{eqnarray*}
Therefore,
$$
(-\Delta)^s u^+ +c(x) u^+ \leq (-\Delta)^s u +c(x) u \leq 0 \mbox{ in } D^+.
$$

 This verifies (\ref{eq2.1}).
\medskip

\emph{Step 2.}
\medskip

 Step 1 shows that $u^+$ is sub-harmonic in $D^+$. In this step, we compare $u^+$ with $s$-harmonic function $\hat{u}.$ We show that
\begin{eqnarray}\label{eq2.3}
 u^+(x)\leq \hat{u}(x) \;\mbox{ in }\; B_r(0)\; \mbox{ for all } \;r>0,
\end{eqnarray}
 where
 \begin{eqnarray}\label{hat}
   \hat{u}(x)=\left\{\begin{array}{ll}
\int_{|y|>r}P_r(y, x)u^+(y) dy, &|x|<r, \\
 u^+(x),& |x|\geq r,
\end{array} \right.
 \end{eqnarray}
Here, $P_r(y, x)$ is the so-called Poisson kernel defined by
\begin{eqnarray*}
P_r(y, x)=\left\{\begin{array}{ll}
B(n, s)\left(\frac{r^2-|x|^2}{|y|^2-r^2}\right)^s \frac{1}{|x-y|^n}, &|y|>r, \\
0,& |y|< r.
\end{array} \right.
\end{eqnarray*}
with $B(n,s)= \frac{\Gamma(n/2)}{\pi^{\frac{n}{2}+1}}\sin(\pi s).$

It is known that (see \cite{CLM}) $\hat{u}(x)$ satisfies
\begin{eqnarray*}
(-\lap)^s\hat{u}(x) =0, \;\;\;  |x|<r.
\end{eqnarray*}

Set
\begin{eqnarray*}v(x)=u^+(x)-\hat{u}(x),
\end{eqnarray*}
obviously, $v(x)\equiv 0$ in $B_r^c(0).$ If (\ref{eq2.3}) does not hold, then there exists a point $\bar{x} \in B_r(0)$ such that
\begin{eqnarray*}
v(\bar{x})=\sup_{B_r(0)} v(x) >0.
\end{eqnarray*}
If $\bar {x}  \in (D^+)^c,$ we have
\begin{eqnarray*}v(\bar{x})=u^+(\bar{x}) -\hat{u}(\bar{x})=-\hat{u}(\bar{x})\leq 0.\end{eqnarray*}
 Therefore,
$\bar{x } \in D^+$.

By (\ref{eq2.1}), we have
\begin{eqnarray*}
(-\Delta)^s v(\bar{x}) = (-\Delta)^s u^+(\bar{x})-(-\Delta)^s \hat{u} (\bar{x})\leq 0.
\end{eqnarray*}
Since $\bar{x}$ is the maximum point of $v$ in $B_r(0),$ by a direct calculation, we have
\begin{eqnarray*}
(-\lap)^s v(\bar{x}) = C_{n,s} PV \int_{\mathbb{R}^n}\frac{v(\bar{x}) -v(y)}{|\bar{x}-y|^{n+2s}}dy>0.
\end{eqnarray*}
\medskip
It is a contradiction. Hence, we arrived at (\ref{eq2.3}).

\emph{Step 3.}
\medskip

 Based on the properties of $u^+$  derived in the previous two steps, we now show that $$u^+(x)\equiv 0 \;\mbox{ in } \;D.$$

 Suppose not,  since $u$ is bounded from above, we have
\begin{eqnarray}\label{eq2.2}
+\infty> A:=\sup_{{\mathbb{R}^n}} u^+>0.
\end{eqnarray}
Then for any $0<\gamma <1$, there exists a point $x_0 \in D$ such that
\begin{eqnarray}\label{x0}
u^+ (x_0)\geq \gamma A.
\end{eqnarray}
For convenience, by a translation we may assume that $x_0=0 \in D,$ we only need to replace $x-x_0$ by $x.$
We will show that $u^+(0)<(1-C)A$ for a positive constant $C$ independent of $\gamma$  to derive a contradiction.
\medskip

Define
\begin{eqnarray*}
{\cal E} _{2s}^{(r)}(x)=\left\{\begin{array}{ll}
0, &|x|<r, \\
 B(n,s)\frac{r^{2s}}{(|x|^2-r^2)^s |x|^n},& |x|>r,
\end{array} \right.
\end{eqnarray*}
then
\begin{eqnarray*}
({\cal E} _{2s}^{(r)} \ast u^+) (x) =B(n,s) \int_{|y-x|>r}\frac{r^{2s}}{(|y-x|^2-r^2)^s|x-y|^n}u^+(y)dy,
\end{eqnarray*}
where $\ast$ denotes the convolution.
Recalling the definition of $\hat{u}$ in (\ref{hat}), we have
\begin{eqnarray*}
\hat{u}(x)=B(n,s) \int_{|y|>r}\frac{(r^2-|x|^2)^s}{(|y|^2-r^2)^s|x-y|^n}u^+(y)dy,\;\; |x|<r.
\end{eqnarray*}
Then \begin{eqnarray*}
\hat{u}(0)=({\cal E}_{2s}^{(r)} \ast u^+) (0).
\end{eqnarray*}

 By (\ref{eq2.3}), we have
\begin{eqnarray*}
u^+(0)\leq \hat{u}(0).
\end{eqnarray*}
Therefore,
\begin{eqnarray}\label{hat1}
 u^+(0)\leq ({\cal E} _{2s}^{(r)} \ast u^+) (0).
\end{eqnarray}
By direct calculations, we know that
\begin{eqnarray}\label{hat2}
\int_{|x|>r}{\cal E} _{2s}^{(r)}(x)dx=B(n,s) \int_{|x|>r}\frac{r^{2s}}{(|x|^2-r^2)^s|x|^n}dy=1.
\end{eqnarray}
Since
$$u^+(y)\leq A, \;\; x\in D \mbox{~and~} u^+(y)=0, \;\; x \in D^c,$$
we obtain from  (\ref{domain}) that there exists $k_0>0$ such that
\begin{eqnarray}\label{hat3}
\mid D^c\cap(B_{2^{i+1}r}\backslash B_{2^{i}r})\mid  \geq \frac{C_0}{2}\mid B_{2^{i+1}r}\backslash B_{2^{i}r}\mid, \;\; \mbox{for any } \; i\geq k_0.
\end{eqnarray}
 Combining (\ref{hat1})-(\ref{hat3}), we derive
\begin{eqnarray}\label{model}
u^+(0) & \leq & B(n,s) \int_{|y|>r}\frac{r^{2s}}{(|y|^2-r^2)^s|y|^n}u^+(y)dy\nonumber\\
&=&B(n,s) \int_{D\cap\{|y|>r\}}\frac{r^{2s}}{(|y|^2-r^2)^s|y|^n}u^+(y)dy\nonumber\\
&&+AB(n,s) \int_{D^c\cap\{|y|>r\}}\frac{r^{2s}}{(|y|^2-r^2)^s|y|^n}dy\nonumber\\
&&-AB(n,s) \int_{D^c\cap\{|y|>r\}}\frac{r^{2s}}{(|y|^2-r^2)^s|y|^n}dy\nonumber\\
&\leq& AB(n,s) \int_{\{|y|>r\}}\frac{r^{2s}}{(|y|^2-r^2)^s|y|^n}dy\nonumber\\
&&-AB(n,s) \int_{D^c\cap\{|y|>r\}}\frac{r^{2s}}{(|y|^2-r^2)^s|y|^n}dy\nonumber\\
&=&A-AB(n,s) \int_{D^c\cap\{|y|>r\}}\frac{r^{2s}}{(|y|^2-r^2)^s|y|^n}dy\nonumber\\
&\leq& A-AB(n,s)\sum\limits_{i = k_0}^{\infty}\int_{D^c\cap({B_{2^{i+1}r}}\backslash B_{2^ir})}\frac{r^{2s}}{(|y|^2-r^2)^s|y|^n}dy\nonumber\\
  &\leq & A-CA\sum\limits_{i = k_0}^{\infty}\frac{r^{2s}}{(4^{i+1}r^2-r^2)^s2^{n(i+1)}r^n}  \mid D^c\cap(B_{2^{i+1}r}\backslash B_{2^{i}r})\mid \nonumber\\
  &\leq&  A-CA\sum\limits_{i = k_0}^{\infty}\frac{1}{(4^{i+1}-1)^s2^{n(i+1)}}  2^{in}\nonumber\\
   &=&(1-C_1)A,
\end{eqnarray}
with $C_1$ independent of $\gamma.$ Taking $\gamma$ closly to $1$ in (\ref{x0}), we derive a contradiction. Therefore we must have
\begin{eqnarray*}
u^+(x)=0, \;\;  x \in D,
\end{eqnarray*}
and therefore
\begin{eqnarray*}
u\leq 0, \;\; x \in D.
\end{eqnarray*}

This completes the proof of Theorem \ref{mthm3}.
\medskip

\subsection{Some immediate applications of the maximum principle}

It is well-known that {\em{ maximum principles}} play important roles in the analysis of the corresponding PDEs. Here we list two simple examples in the following two lemmas.

Define $$\mathbb{R}_+^n=\{x=(x_1, \cdots , x_n) \mid x_n>0\}.$$

\begin{lem}\label{A}
Let $u \in {\cal L}_{2s}(\mathbb{R}^n)\cap C_{loc}^{1,1}( \mathbb{R}^n_+)\cap C( \overline{{\mathbb{R}^n_+}})$ be a bounded non-negative solution of
\begin{eqnarray*}
\left\{\begin{array}{ll}
(-\Delta)^s u =f(u),&x \in \mathbb{R}^n_+;\\
u>0,&x \in \mathbb{R}^n_+;\\
u= 0, & x \notin  \mathbb{R}^n_+.
\end{array}
\right.
\end{eqnarray*}
The given function $f=f(u)$ is non-increasing.
Then $u$ is strictly monotone increasing in the $x_n$ direction.
\end{lem}
\medskip

{\bf{Proof.}}
Let $$D=\mathbb{R}_+^n,$$ $$u^\tau(x)=u(x', x_n+\tau)$$ and $$U^\tau(x)=u(x)-u^\tau(x).$$
Since $f(u)$ is non-increasing,
 $$f(u(x))-f(u^\tau(x))\leq 0 \mbox{ in } D \mbox{ where } U^\tau(x)>0.$$
Therefore,
$$
(-\lap)^s U^\tau(x)\leq 0 \mbox{ in } D \mbox{ where } U^\tau(x)>0.
$$

In addition,
 $$U^\tau(x)\leq 0, \;\; x \notin \mathbb{R}^n_+.$$
 It follows from Theorem \ref{mthm3} that
 $$U^\tau(x) \leq 0, \;\; x \in \mathbb{R}^n_+,\;\; \mbox{for any }\tau >0.$$

Now we show that
 $$U^\tau(x) < 0, \;\; x \in \mathbb{R}^n_+,\;\; \mbox{for any }\tau >0.$$
If not, there exists a point $x^0 \in \mathbb{R}^n_+$ such that
$$
U^\tau(x^0) =0=\max_{\mathbb{R}^n}U^\tau(x).
$$
Since $U^\tau(x)\not\equiv 0$ in $\mathbb{R}^n,$ we have
$$
(-\Delta)^sU^\tau(x^0) =C_{n, s} PV \int_{\mathbb{R}^n}\frac{0-U^\tau(y)}{|x^0-y|^{n+2s}}dy<0,
$$
while it contradicts
$$
(-\Delta)^sU^\tau(x^0) =f(u(x^0))-f(u^\tau(x^0))=0.
$$
Therefore $u$ is strictly monotone increasing in the $x_n$ direction.

This proves Lemma \ref{A}.
\medskip

\begin{lem}\label{A2}
Suppose that $u \in {\cal L}_{2s}(\mathbb{R}^n)\cap C_{loc}^{1,1}(\mathbb{R}^n)$ is a solution of
$$
(-\lap)^s u(x) = f(u(x)), \;\;\;\;  x \in \mathbb{R}^n,
$$
and
$$|u(x)| \leq 1, \;\;\; \forall \, x \in \mathbb{R}^n, $$
\be \lim_{x_n \ra \pm \infty} u(x',x_n) =  \pm 1 \mbox{ uniformly for all } x' \in \mathbb{R}^{n-1}.
\label{1u}
\ee
Assume that there exists $\delta >0$ such that
\begin{eqnarray}\label{1f}
f \;\mbox{ is nonincreasing on }\;[-1, -1+\delta] \;\mbox{ and on } \;[1-\delta, 1].
\end{eqnarray}
Then $u^\tau(x)\geq u(x)$ in $\mathbb{R}^n$ for sufficiently large $\tau$.
\end{lem}
\medskip

We prove this Lemma by using Theorem \ref{mthm3} {\em{(the maximum principle)}}.
\medskip

\textbf{Proof.} Denote
$$
U^\tau(x)=u(x)-u^\tau(x),\;\; x\in \mathbb{R}^n.
$$
We prove that for sufficiently large $\tau,$
\begin{eqnarray}\label{LA}
U^\tau(x)\leq 0 \mbox{ in } \mathbb{R}^n.
\end{eqnarray}
If not, then
\begin{eqnarray}\label{CA}
\sup_{\mathbb{R}^n}U^\tau(x):=A>0.
\end{eqnarray}

Consider the function $U^\tau(x)- \frac{A}{2}.$ First,
we obtain from (\ref{1u}) that there exists a constant $a>0$ such that
$$u(x', x_n) \geq 1-\delta, \; \mbox{ for } x_n \geq a$$
 and
  $$u(x', x_n) \leq -1+\delta, \; \mbox{ for } x_n \leq -a,$$
and  there exists a positive constant $M>a$ such that
\begin{eqnarray}\label{IA}
U^\tau(x)-\frac{A}{2}\leq 0, \;\; x \in \mathbb{R}^{n-1}\times [M, +\infty),
\end{eqnarray}

Second, for any $\tau \geq 2a,$ no matter where $x$ is, one of the points $x$ and $x+(0', \tau)$ is in the domain $\{x: |x_n|\geq a\}$. We have either
\begin{eqnarray}\label{le1}
u^\tau (x', x_n) \geq 1-\delta \;\;  ( \mbox{ if } x_n \geq -a),
\end{eqnarray}
or
\begin{eqnarray}\label{le2}
u (x', x_n) \leq -1+\delta \;\; ( \mbox{ if } x_n \leq -a),
\end{eqnarray}
 where $f$ is non-increasing. Therefore,
$$
f(u(x))-f(u^\tau(x))\leq 0 \mbox{ at the point in } D=\mathbb{R}^{n-1}\times (-\infty, M) \mbox{ where } u(x)>u^\tau(x).
$$
Consequently,
$$
(-\Delta)^s\left(U^\tau-\frac{A}{2}\right)\leq 0 \mbox{ at the point in } D \mbox{ where } U^\tau(x)-\frac{A}{2}>0.
$$
It follows from this, (\ref{IA}) and Theorem  \ref{mthm3} that

 $$U^\tau(x)- \frac{A}{2}\leq 0, \;\; x\in \mathbb{R}^n.$$
This contradicts (\ref{CA}), thus we derive (\ref{LA}) and complete the proof of Lemma \ref{A2}.

\section{The monotonicity of solutions in $\mathbb{R}^n$}

In this section, we prove Theorem \ref{mthm4}.
\medskip

{\em Outline of the proof.}
\smallskip

Write $x=(x', x_n)$. For any $\tau \in \mathbb{R}$, define $$u^\tau (x)= u(x', x_n+\tau),$$
and
$$U^\tau (x)=u(x)-u^\tau (x).$$
Our aim is to prove that for any  $\tau>0$, we have
\begin{eqnarray*}
U^\tau (x)< 0 \mbox{ in } \mathbb{R}^n .
\end{eqnarray*}
This can be immediately deduced from
\begin{eqnarray}\label{FA}
U^\tau (x)\leq 0 \mbox{ in } \mathbb{R}^n
\end{eqnarray}
by a similar argument as in proving strong maximum principles.

First, we derive (\ref{FA})  for sufficiently large $\tau$ as a consequence of Lemma \ref{A2}. This provides a starting point for sliding the domain.

Then we decrease $\tau$ to its limit as long as inequality (\ref{FA}) holds. Define
$$\tau_0= \inf \, \{\tau \mid U^\tau(x) \leq 0, x \in \mathbb{R}^n\}.$$
We show that
 $$ \tau_0=0.$$
 Otherwise, we prove that $\tau_0$ can be decreased a little bit while inequality
(\ref{FA}) is still valid, which would contradict the definition of $\tau_0.$
From the continuity of $U^\tau$ with respect to $\tau,$ this can be realized if we have
\begin{eqnarray}\label{supa}
\sup_{\mathbb{R}^n} U^{\tau_0}(x)<0.
\end{eqnarray}
Unfortunately, this is not true since
$$
\lim_{x_n\rightarrow \pm\infty} U^{\tau_0}(x)=0.
$$
Hence, instead of (\ref{supa}), we prove
\begin{eqnarray}\label{supb}
\sup_{|x_n|\leq a } U^{\tau_0}(x)<0.
\end{eqnarray}
Since the maximum principle (Theorem \ref{mthm3}) is not applicable here, we introduce a generalize average inequality to derive a contradiction at the maximum points of a sequence of auxiliary functions.

(\ref{supb}) implies immediately that for some small $\varepsilon >0,$
$$
U^\tau(x)\leq 0, \;  |x_n|\leq a,\; \tau \in (\tau_0-\varepsilon, \tau_0].
$$
Then what left is to show that
$$
U^\tau(x)\leq 0, \;  |x_n|> a,\; \tau \in (\tau_0-\varepsilon, \tau_0].
$$
This is similar to the proof of Lemma \ref{A2}.

\medskip

{\em Now we carry out the details of the proof.}
\smallskip

We divide the proof into three steps.
\medskip

\emph{Step 1.}
\medskip

As a consequence of Lemma \ref{A2}, we have
\begin{eqnarray}\label{pstart}
U^\tau (x) \leq 0,  \;\mbox{ in } \;\mathbb{R}^n,\; \tau \geq 2a.
\end{eqnarray}

\medskip

\emph{Step 2.}
\medskip

 (\ref{pstart}) provides a starting point, from which we can carry out the sliding. From $\tau = 2a$, we decrease $\tau,$ and show that for any $0<\tau<2a$, we also have
\begin{eqnarray}\label{Pdecrease}
U^\tau(x)\leq 0,\;\; x \in \mathbb{R}^n.
\end{eqnarray}

Define
$$\tau_0= \inf \, \{\tau \mid U^\tau(x) \leq 0,\; x \in \mathbb{R}^n\}.$$

We prove that  $\tau_0=0$. Otherwise,
 we show that $\tau_0$ can be decreased a little bit while inequality
(\ref{Pdecrease}) is still valid.
\medskip

(i) We first prove that
\begin{eqnarray}\label{tau0}
\sup_{\mathbb{R}^{n-1}\times[-a, a]}U^{\tau_0}(x)<0.
\end{eqnarray}

If not, then
 $$\mathop{\sup}\limits_{\mathbb{R}^{n-1}\times[-a, a]}U^{\tau_0}(x)=0,$$
 and there exists a sequence
 $$\{x^k\} \subset \mathbb{R}^{n-1}\times[-a, a],\; k=1, 2, \cdots,$$
 such that $$U^{\tau_0}(x^k) \rightarrow 0, \mbox{ as } k\rightarrow \infty.$$

Let
\begin{eqnarray}\label{eta}
 \eta(x)=\left\{\begin{array}{ll}
a e^{\frac{1}{|x|^2-1}}, &|x|<1, \\
  0,& |x|\geq 1,
\end{array} \right.
\end{eqnarray}
taking $a=e$ such that $\eta (0)=\mathop{\max}\limits_{\mathbb{R}^n}\eta (x)=1.$

Set
$$
\psi_k(x)=\eta (x-x^k).
$$
 Then there exists a sequence  $\{\varepsilon_k\}\rightarrow 0$ such that
 $$
 U^{\tau_0}(x^k)+\varepsilon_k \psi_k(x^k)> 0.
 $$
For any $x \in \mathbb{R}^n \backslash B_{1}(x^k),$ notice that $U^{\tau_0}(x)\leq 0$ and $ \psi_k(x)=0,$
we have
$$
U^{\tau_0}(x^k) +\varepsilon_k \psi_k(x^k) > U^{\tau_0}(x) +\varepsilon_k \psi_k(x),\; \mbox{ for any }\; x \in \mathbb{R}^n \backslash  B_1(x^k).
$$
It follows that there exists a point $\bar{x}^k \in B_1(x^k)$ such that
\begin{eqnarray}\label{Pwmaximum}
 U^{\tau_0}(\bar{x}^k) +\varepsilon_k \psi_k(\bar{x}^k)=\max_{\mathbb{R}^n}(  U^{\tau_0}(x) +\varepsilon_k \psi_k(x))>0.
\end{eqnarray}

In addition, it can be seen from
$$
U^{\tau_0}(\bar{x}^k) +\varepsilon_k \psi_k(\bar{x}^k)\geq U^{\tau_0}(x^k) +\varepsilon_k \psi_k(x^k),
$$
and $\psi_k(\bar{x}^k)\leq \psi_k(x^k)$ that
$$
0\geq U^{\tau_0}(\bar{x}^k) \geq U^{\tau_0}(x^k).
$$
Therefore,
\begin{eqnarray}\label{F0}
U^{\tau_0}(\bar{x}^k) \rightarrow 0, k\rightarrow\infty.
\end{eqnarray}

It follows from the continuity of $f$ that
 \begin{eqnarray}\label{Pupper}
(-\Delta)^s (U^{\tau_0}+\varepsilon_k \psi_k)(\bar{x}^k)=f(u(\bar{x}^k))-f(u^{\tau_0}(\bar{x}^k)) +C \varepsilon_k \rightarrow 0,\; k\rightarrow\infty.
\end{eqnarray}
In order to apply the maximum principle (Theorem \ref{mthm3}) which is based on the sub-average inequality (\ref{hat1}). We require
$$
(-\Delta)^s(U^{\tau_0}+\varepsilon_k \psi_k)^+(x)\leq 0 \mbox{ in } \mathbb{R}^n.
$$
This is difficult to realize in our situation. For this reason, we introduce the following more general inequality.

\begin{lem}[A generalize average inequality]\label{GAM}
Suppose that  $u \in {{\cal{L}}_{2s}(\mathbb{R}^n)} \cap C_{loc}^{1,1}(\mathbb{R}^n ),$ and $\bar{x}$ is a maximum point of $u$ in $\mathbb{R}^n.$ Then for any $r>0,$ we have
\begin{eqnarray}\label{IGA}
\frac{C_0}{C_{n, s}}r^{2s}(-\Delta)^{s}u(\bar{x})+ C_0 \int_{B_r^c(\bar{x})}\frac{r^{2s}}{|\bar{x}-y|^{n+2s}}u(y)dy\geq u(\bar{x}),
\end{eqnarray}
where $C_0$ satisfies
$$
 C_0 \int_{B_r^c(\bar{x})}\frac{r^{2s}}{|\bar{x}-y|^{n+2s}}dy=1.
$$
\end{lem}

{\bf{Proof.}} By the definition, we have
\begin{eqnarray*}
&&(-\Delta)^{s}u(\bar{x})\\
&=&C_{n, s}PV \int_{\mathbb{R}^n}\frac{u(\bar{x})-u(y)}{|\bar{x}-y|^{n+2s}}dy\\
&\geq& C_{n, s} \int_{B_r^c(\bar{x})}\frac{u(\bar{x})-u(y)}{|\bar{x}-y|^{n+2s}}dy\\
&=& -C_{n, s}\int_{B_r^c(\bar{x})}\frac{u(y)}{|\bar{x}-y|^{n+2s}}dy+ C_{n, s}u(\bar{x}) \int_{B_r^c(\bar{x})}\frac{1}{|\bar{x}-y|^{n+2s}}dy\\
&=& -C_{n, s}\int_{B_r^c(\bar{x})}\frac{u(y)}{|\bar{x}-y|^{n+2s}}dy+ C_{n, s} \int_{B_1^c(0)}\frac{1}{|y|^{n+2s}}dy\frac{u(\bar{x})}{r^{2s}}.\\
\end{eqnarray*}
Multiply both sides by $\frac{C_0 }{ C_{n, s}}r^{2s},$ where
$$C_0=\frac{1}{\int_{B_1^c(0)}\frac{1}{|y|^{n+2s}}dy},
$$
we derive (\ref{IGA}).

This completes the proof of Lemma \ref{GAM}.
\medskip

\begin{mrem} Recalling the definition of ${\cal E} _{2s}^{(r)} \ast u^+ $ in (\ref{hat1}),
\begin{eqnarray*}
({\cal E} _{2s}^{(r)} \ast u^+) (x) =B(n,s) \int_{B_r(x)}\frac{r^{2s}}{(|y-x|^2-r^2)^s|x-y|^n}u^+(y)dy,
\end{eqnarray*}
The proof of maximum principle (Theorem \ref{mthm3}) is adopted by the idea of integral average, where the integral kernel is $\frac{r^{2s}}{(|y-x|^2-r^2)^s|x-y|^n}.$ Similarly, we use the idea of integral average in this case, the corresponding integral kernel is $\frac{r^{2s}}{|\bar{x}-y|^{n+2s}}$. Both of them are equivalent to
 $\frac{r^{2s}}{|y|^{n+2s}}$ if $|y|$ is  sufficiently large.
\end{mrem}

Now we continue the proof of Theorem \ref{mthm4}.

For the function $U^{\tau_0}+\varepsilon_k \psi_k,$ using Lemma \ref{GAM}, we derive
\begin{eqnarray*}
C_1(-\Delta)^{s}(U^{\tau_0}+\varepsilon_k \psi_k)(\bar{x}^k)+ C_2 \int_{B_2^c(\bar{x}^k)}\frac{(U^{\tau_0}+\varepsilon_k \psi_k)(y)}{|\bar{x}^k-y|^{n+2s}}dy\geq (U^{\tau_0}+\varepsilon_k \psi_k)(\bar{x}^k).
\end{eqnarray*}

Combining this with  (\ref{F0}) and (\ref{Pupper}), and the fact that $\varepsilon_k\rightarrow 0,$ we arrive at
\begin{eqnarray}\label{Plower}
0 \leftarrow  \int_{B_2^c(\bar{x}^k)} \frac{U^{\tau_0} (y) }{|\bar{x}^k -y|^{n+2s}}dy
= \int_{B_2^c(0)} \frac{U^{\tau_0} (z+\bar{x}^k) }{|z|^{n+2s}}dz ,\mbox{ as } k\rightarrow\infty.
\end{eqnarray}

Denote $$u_k(x)=u(x+\bar{x}^k) \mbox{ and } U^{\tau_0}_k(x)=U^{\tau_0}(x+\bar{x}^k).$$
Since $u(x)$ is uniformly continuous, by Arzel\`{a}-Ascoli theorem, up to extraction of a subsequence, we have
$$
u_k(x)\rightarrow u_\infty (x) \mbox { in } \mathbb{R}^n, \mbox{ as } k\rightarrow \infty,
$$
and by (\ref{Plower}),
$$
U^{\tau_0}_k(x)\rightarrow  0, \;\; x \in B_2^c(0), \mbox{ uniformly, as } k\rightarrow\infty.
$$
Therefore,
$$
 u_\infty (x)-u^{\tau_0}_\infty (x)\equiv 0,\;\; x \in B_2^c(0).
$$

Since $\{\bar{x}^k_n\}$ is bounded, we obtain from (\ref{eq1.2}) that
\begin{eqnarray}\label{i3}
u_\infty(x', x_n)  \mathop {\longrightarrow }\limits_{{x_n\rightarrow \pm \infty}}\pm 1 \;\mbox{ uniformly in }\; x'=(x_1, \cdots, x_{n-1}).
\end{eqnarray}
Therefore,
\begin{eqnarray}\label{i4}
u_\infty(x', x_n)&=&u_\infty(x', x_n+\tau_0)=u_\infty(x', x_n+2\tau_0)\nonumber\\
&=&\cdots=u_\infty(x', x_n+k\tau_0),
\end{eqnarray}
for any $k \in\mathbb{N}.$

From (\ref{i3}) one can take $x_n$ sufficiently negative to make $u_\infty(x', x_n)$ close to $-1,$ and then take $k$ sufficiently large to make $u_\infty(x', x_n+k\tau_0)$ close to $1$, this is impossible.
Hence (\ref{tau0}) must hold.
\medskip

(ii) We prove that, there exists an $\varepsilon >0$, such that
\begin{eqnarray}\label{tau minus}
U^{\tau}(x)\leq0, \;\; \forall x \in \mathbb{R}^n, \;\; \forall \tau \in (\tau_0-\varepsilon,\tau_0].
\end{eqnarray}

First, (\ref{tau0}) implies immediately  that there exists a small $\varepsilon >0$ such that
\begin{eqnarray}\label{ib}
\sup_{\mathbb{R}^{n-1}\times[-a, a]}U^{\tau}(x)<0,\;\; \forall \tau \in (\tau_0-\varepsilon, \tau_0].
\end{eqnarray}
Therefore, we only need to prove that
\begin{eqnarray}\label{NFB}
\sup_{\mathbb{R}^{n}\backslash (\mathbb{R}^{n-1}\times[-a, a])}U^{\tau}(x)\leq 0,\;\; \forall \tau \in (\tau_0-\varepsilon, \tau_0].
\end{eqnarray}
If not, then
\begin{eqnarray}\label{onestar}
\sup_{\mathbb{R}^{n}\backslash (\mathbb{R}^{n-1}\times[-a, a])}U^{\tau}(x):=A> 0,\;\; \forall \tau \in (\tau_0-\varepsilon, \tau_0].
\end{eqnarray}

By the asymptotic condition (\ref{eq1.2}), we may assume that there exists $M>a$ such that
\begin{eqnarray}\label{NF0}
U^{\tau}(x)\leq \frac{A}{2}, \;\; x \in \mathbb{R}^n \backslash (\mathbb{R}^{n-1}\times [-M, M]),\;\; \forall \tau \in (\tau_0-\varepsilon, \tau_0].
\end{eqnarray}

Denote
 $$D_1=(\mathbb{R}^{n-1}\times (a, M)) \cup (\mathbb{R}^{n-1}\times (-M, -a)),$$
 then
$$
D_1^c=\mathbb{R}^{n}\backslash (\mathbb{R}^{n-1}\times[-M, M])\cup (\mathbb{R}^{n-1}\times [-a, a]).
$$
Consider the function $U^{\tau}(x)- \frac{A}{2}.$ First, by (\ref{ib}) and (\ref{NF0}), we have
\begin{eqnarray}\label{EU}
U^{\tau}(x)- \frac{A}{2}\leq 0, \;  x \in D_1^c,\;\; \forall \tau \in (\tau_0-\varepsilon, \tau_0].
\end{eqnarray}

Second, for any $\tau \in (\tau_0-\varepsilon, \tau_0],$ at the points in $D_1^c$ where $U^{\tau}(x)- \frac{A}{2}>0,$ we have $u(x)\geq u^{\tau}(x).$ If $x \in \mathbb{R}^{n-1}\times (a, M),$ we have $u(x), u^{\tau}(x)\geq 1-\delta,$
and $f(u(x))\leq f(u^{\tau}(x))$  due to the monotonicity of $f.$
If $x \in \mathbb{R}^{n-1}\times (-M, -a),$ we have $-1+\delta\geq u(x)\geq u^{\tau}(x).$ Then  we also have $f(u(x))\leq f(u^{\tau}(x))$ due to the monotonicity of $f.$

 Therefore,
\begin{eqnarray}\label{IU}
(-\Delta)^s \left(U^{\tau}(x)- \frac{A}{2}\right)=f(u(x))- f(u^{\tau}(x))\leq 0,\;\mbox{in } \left\{x\in D_1^c \mid U^{\tau}(x)> \frac{A}{2}\right\}.
\end{eqnarray}
Combining (\ref{EU}) with (\ref{IU}), and by Theorem \ref{mthm3}, we arrive at
$$
U^{\tau}(x)- \frac{A}{2}\leq 0, \;  x \in \mathbb{R}^n,\;\; \forall \tau \in (\tau_0-\varepsilon, \tau_0].
$$
This contradicts (\ref{onestar}). Therefore, (\ref{NFB}) is correct.
This proves (\ref{tau minus}) which contradicts the definition of $\tau_0$. Hence we obtain (\ref{Pdecrease}).
\medskip

\emph{Step 3.}
\medskip

In this step, we show that $u$ is strictly increasing with respect to $x_n$ and $u(x)$ depends on $x_n$ only.
\medskip

 Combining these two steps above, we have derived that
 \begin{eqnarray*}
U^\tau (x) \leq 0  \;\mbox{ in } \;\mathbb{R}^n, \;\;  \forall\tau >0.
\end{eqnarray*}
Now based on this, if ${U^\tau}(x^o) = 0 $ at some point $x^o \in \mathbb{R}^n,$ then $x^o$ is a maximum point of $U^\tau$ in $\mathbb{R}^n,$ and
$$
(-\Delta)^s {U^\tau}(x^o)=f (u(x^o))-f ({u^\tau}(x^o))=0.
$$
Since $U^\tau(y) \not\equiv 0$ in $\mathbb{R}^n$, by a direct calculation, we have
\begin{eqnarray*}
&&(-\Delta)^s {U^\tau}(x^o)\\
&=&C_{n, s} PV \int_{\mathbb{R}^n}\frac{{U^\tau}(x^o)-{U^\tau}(y)}{|x^o-y|^{n+2s}}dy\\
&= &C_{n, s} PV \int_{\mathbb{R}^n}\frac{-U^{\tau}(y)}{|x^o-y|^{n+2s}}dy\\
&>&0.
\end{eqnarray*}
It is a contradiction. Therefore, we have
 \begin{eqnarray*}
U^\tau (x) < 0  \;\mbox{ in } \;\mathbb{R}^n, \;\; \mbox{ for }  \forall\tau >0.
\end{eqnarray*}
This implies that $u$ is strictly increasing with respect to $x_n.$
 \smallskip

Now we claim that $u(x)$ depends on $x_n$ only.

In fact, it can be seen from the above process that the argument still holds if we replace $u^\tau (x)$ by $u(x+\tau \nu)$, where $\nu=(\nu_1,\cdots,\nu_n)$ with $\nu_n >0$ is an arbitrary vector pointing upward.
Applying the similar argument as in Step 1 and 2, we can derive that, for each of such $\nu$,
$$u(x+\tau \nu) > u(x), \;\; \forall \, \tau > 0, \, x \in \mathbb{R}^n.$$

Let $\nu_n \ra 0$,  from the continuity of $u,$ we deduce that $$u(x+\tau \nu)\geq u(x)$$ for arbitrary $\nu$ with $\nu_n=0$. By replacing $\nu$ by $-\nu,$ we find that $$u(x+\tau \nu) = u(x)$$ for arbitrary $\nu$ with $\nu_n=0$, this means that $u$ is independent of $x'$, hence $u(x)=u(x_n)$.

This completes the proof of Theorem \ref{mthm4}.

\bigskip

{\em Authors' Addresses and E-mails:}
\medskip

Leyun Wu

School of Mathematical Sciences, MOE-LSC,

Shanghai Jiao Tong University

Shanghai, China and

Department of Mathematical Sciences

Yeshiva University

New York, NY, 10033 USA

leyunwu@sjtu.edu.cn
\medskip

Wenxiong Chen

Department of Mathematical Sciences

Yeshiva University

New York, NY, 10033 USA

wchen@yu.edu
\medskip

\end{document}